\documentclass[12pt]{article}

\usepackage{theorem,amssymb}

\advance \topmargin by -\headheight
\advance \topmargin by -\headsep
\evensidemargin \oddsidemargin
\author{J.-P. Allouche \\
CNRS, LRI \\
B\^atiment 490 \\
F-91405 Orsay Cedex (France) \\
{\tt allouche@lri.fr} \\
\and
G. Skordev \\
CEVIS, Universit\"at Bremen \\
Universit\"atsallee 29 \\
D-28359 Bremen (Germany) \\
{\tt skordev@cevis.uni-bremen.de} \\
}

\title{Von Koch and Thue-Morse revisited}

\date{ }

\def \proof{\bigbreak\noindent{\it Proof.\ \ }}

\def \endpf{{\ \ $\Box$ \medbreak}}

\theorembodyfont{\rm}

\newtheorem{definition}{Definition}

\begin{document}

\maketitle

\begin{abstract}
We revisit the relation between the von Koch curve and the Thue-Morse 
sequence given in a recent paper of Ma and Goldener by relating their 
study to papers written by Coquet and Dekking at the beginning of 
the 80s. We also emphasize that more general links between fractal 
objects and automatic sequences can be found in the literature.
\end{abstract}

\section{Introduction}

The recent paper of Ma and Holdener \cite{MH} gives an interesting explicit
relation between the von Koch curve and the Thue-Morse sequence. We show here 
that such a link was already noted in two papers dated 1982--1983, one by 
Coquet \cite{Coquet} (where the name ``von Koch'' is not explicitly mentioned, 
but a construction \`a la von Koch is given \cite[p.~111]{Coquet}), the other 
by Dekking, see \cite[p.~32-05 and p.~32-06]{Dekking-bx}.
We explain the similarity of the points of view of these three papers through 
the approach of Dekking in \cite{Dekking1} (where the Thue-Morse sequence
does not appear). Actually the fundamental object is the sequence
$((-1)^{s_2(n)} j^n)_{n \geq 0}$ and its summatory function, where $s_2(n)$ 
is the sum of the binary digits of the integer $n$, hence 
$((-1)^{s_2(n)})_{n \geq 0}$ is the $\pm 1$ Thue-Morse sequence, and
$j = e^{2i\pi/3}$. The sum $\sum_{n < N} (-1)^{s_2(n)} j^n$ was explicitly 
introduced by Coquet in \cite{Coquet} to compute $\sum_{n < N} (-1)^{s_2(3n)}$, 
and in a more general form by Dekking in \cite{Dekking-bx}. 
Joining the points of the complex plane that correspond to the sums 
$\sum_{n < N} (-1)^{s_2(n)} j^n$ gives a non-renormalized version of the von 
Koch curve. We end this paper with a quick survey of some general relations 
between automatic sequences and fractal objects.

\bigskip

We suppose that the reader is familiar with the construction of the von Koch 
curve (see in particular the original paper of von Koch \cite{Koch} or the 
book of Mandelbrot \cite{Mandelbrot}) and with the Thue-Morse sequence (see 
for example \cite{AS1}). We will nevertheless recall the definition of the 
Thue-Morse sequence as well as the definition of morphisms of the free monoid 
and of automatic sequences: for more details the reader can look at \cite{AS2} 
for example.

\section{A few definitions from combinatorics on words}

In what follows an {\it alphabet\,} is a nonempty finite set. Its elements
are called {\it letters}. A {\it word\,} on the alphabet ${\cal A}$ is a 
finite string of elements of ${\cal A}$. The {\it length\,} of a word is 
the number of letters composing this word. The {\it concatenation\,} of two 
words $v = a_1 a_2 \cdots a_i$ and $w = b_1 b_2 \cdots b_j$, where the $a$'s 
and $b$'s are letters is noted $vw$ and defined by 
$vw := a_1 a_2 \cdots a_i b_1 b_2 \cdots b_j$. The {\it free monoid\,} 
generated by ${\cal A}$ is the set of words on ${\cal A}$, including the 
empty word, equipped with the concatenation of words: it is denoted 
by ${\cal A}^*$.

\begin{definition} Let ${\cal A}$ be an alphabet. 

\medskip

A {\em morphism} $\varphi$ on ${\cal A}$ is a map from ${\cal A}^*$ to 
itself that is a homomorphism for the concatenation. (The morphism $\varphi$ 
is determined by the images $\varphi(a)$ for $a \in {\cal A}$.)

\medskip

A morphism $\varphi$ on ${\cal A}$ is called {\it $d$-uniform\,} for some 
integer $d \geq 1$ (or {\it uniform of length $d$\,} or {\it uniform\,} if 
$d$ is clear from the context) if  the words $\varphi(a)$ have length $d$ 
for all $a \in {\cal A}$.
 
\end{definition}

\bigskip

\noindent
{\bf Example} \ \
Let ${\cal A} := \{0, 1\}$. Define the $2$-uniform $\varphi$ on ${\cal A}$ 
by $\varphi(0) := 01$, $\varphi(1) := 10$. Then, denoting by $\varphi^k$ the
$k$th iterate of $\varphi$, we have:
$$
\varphi(0) = 01, \ \ \ \varphi^2(0) = 0110, \ \ \ \varphi^3(0) = 01101001, \ \ \
\ldots
$$
It is easy to see that the sequence of words $(\varphi^k(0))_k$ simply 
converges to the infinite sequence
$$
0 \ 1 \ 1 \ 0 \ 1 \ 0 \ 0 \ 1 \ 1 \ 0 \ 0 \ 1 \ 0 \ 1 \ 1 \ 0 \ \ldots
$$
called the {\it Thue-Morse} sequence on the alphabet $\{0, 1\}$, and that,
for $n \geq 0$, the $n$th term of this sequence is equal to the parity
of $s_2(n)$, the sum of the binary digits of the integer $n$. Replacing 
$0$ by $a$ and $1$ by $b$ yields the Thue-Morse sequence on the alphabet 
$\{a, b\}$ (called the $\pm 1$ Thue-Morse sequence if $a=1$ and $b=-1$)
$$
a \ b \ b \ a \ b \ a \ a \ b \ b \ a \ a \ b \ a \ b \ b \ a \ \ldots
$$

\begin{definition}
Let ${\cal A}$ be an alphabet and $d$ be a positive integer. An infinite 
sequence $(u_n)_{n \geq 0}$ on ${\cal A}$ is called $d$-automatic if there 
exist an alphabet ${\cal B}$, a letter $b \in {\cal B}$, a $d$-uniform 
morphism $\varphi$ on ${\cal B}$ and a map $f : {\cal B} \to {\cal A}$ 
such that the sequence of words $(\varphi^k(b))_k$ simply converges to 
an infinite sequence $(v_n)_{n \geq 0}$ on ${\cal B}$, and such that for 
all $n \geq 0$ we have $u_n = f(v_n)$.
\end{definition}

\section{The result of Ma and Holdener}

The authors of \cite[p.~202]{MH} describe the picture obtained with a ``turtle 
program'' (see \cite{AdS} for example, see also\cite{HW}) based on the Thue-Morse 
sequence. After scaling, the picture converges to the von Koch snowflake. Up to 
designing a ``simplified'' curve, the authors prove in fact that the classical 
iterative construction scheme for the von Koch curve is precisely given by the 
Thue-Morse sequence in this setting. We begin with a simplified formulation of 
their result.

\noindent
{\bf Theorem (Ma, Holdener \cite{MH})} \
Let $\{F, L\}$ be the turtle alphabet where $F$ is {\tt move one unit forward} 
and $L$ is {\tt rotate by $\pi/3$}. Let $TM_k$ be the prefix of length $2^k$ 
of the Thue-Morse sequence on the alphabet $\{F, L\}$. Let 
$W_n := K_{2n}(TM_{4n})$ be defined as follows: the word $TM_{4n}$ is written 
as the concatenation of words from $\{TM_{2n}, \overline{TM}_{2n}\}$, where
$\overline{TM}_{2n}$ is obtained from $TM_{2n}$ by interchanging the symbols 
$F$ and $L$; this concatenation is interpreted as a word on the turtle alphabet 
$\{F, L\}$ thus defining a polygonal line in the plane denoted by 
$\widetilde{W}_n$. Let $S_n$ be the scaling factor 
defined by $S_n := 2/(3^{n-2}-1)$. Then the sequence of polygonal lines 
$S_{2n} \widetilde{W}_n$ converges to the von Koch curve.

\bigskip

\proof This is essentially Theorem~5.0.14 of \cite{MH}. 
\endpf

\bigskip

One sees that, up to scaling, the polygonal line $\widetilde{W}_{n+1}$ is 
obtained from $\widetilde{W}_n$ in four steps: draw $\widetilde{W}_n$, glue 
a copy of $\widetilde{W}_n$ rotated by $-\pi/3$, glue a copy of 
$\widetilde{W}_n$ rotated by $-\pi/3 + 2\pi/3 = \pi/3$, and glue a last copy 
of $\widetilde{W}_n$ rotated by $-\pi/3 + 2\pi/3 -\pi/3 = 0$. In other words, 
the polygonal line $\widetilde{W}_n$ is transformed into $\widetilde{W}_{n+1}$ 
exactly as follows, up to scaling: each segment in $\widetilde{W}_n$ is 
interpreted as a complex number $z$ in the coordinates where the origin is 
an endpoint of this segment and the $x$-axis is parallel to this segment; 
this complex number is then replaced by the four complex numbers 
$z, -jz, -j^2z, z$ (where $j := e^{2i\pi/3}$); the segments corresponding to 
$z, -jz, -j^2z, z$ in the same coordinates are glued consecutively, i.e., the 
segment we started from is replaced by the segments corresponding to the complex 
numbers $z, z-jz, z-jz-j^2z, z-jz-j^2z +z$.

In other words, still ignoring scaling, the curve which is obtained
is the geometric representation of the sum $\sum_{0 \leq k < 4^n} u_k$
where the sequence $(u_n)_{n \geq 0}$ is the fixed point beginning in $1$
of the morphism defined over the alphabet $\{1, -1, j, -j^2\}$ by
$$
\lambda:
\left\{\begin{array}{rrrrrr}
1    &\to &1    \ &-j   \ &-j^2 \ &1    \\
-1   &\to &-1   \ &j    \ &j^2  \ &-1   \\
j    &\to &j    \ &-j^2 \ &-1   \ &j    \\
-j   &\to &-j   \ &j^2  \ &1    \ &-j   \\
j^2  &\to &j^2  \ &-1   \ &-j   \ &j^2  \\
-j^2 &\to &-j^2 \ &1    \ &j    \ &-j^2 \\
\end{array}
\right.
$$
which easily implies by induction on $n$ that
$$
u_n = (-1)^{s_2(n)} j^n
$$
where $s_2(n)$ is the sum of the binary digits of the integer $n$, and 
where thus the sequence $((-1)^{s_2(n)})_{n \geq 0}$ is the $\pm 1$ 
Thue-Morse sequence. Note that, to recover the usual (bounded) von Koch
curve, the appropriate scaling here is $3^{-n}$.

\bigskip

This yields two consequences:

\begin {itemize}

\item The occurrence of the Thue-Morse sequence in \cite{MH} parallels its 
occurrence in \cite{Coquet}, where Coquet, interested in the behavior of
the summatory function $\sum_{k < n} (-1)^{s_2(3k)}$, introduced the sum
$\sum_{k < n} (-1)^{s_2(k)} j^k$. He obtains (page 111) as a by-product 
the von Koch curve, which he calls ``a fractal object built on a classical 
fractal scheme'' (see also \cite{Dekking-bx}).

\item The morphism above is a simplification of the morphism given by Dekking 
in \cite{Dekking1} to construct the von Koch curve. Namely, Dekking considered
the morphism defined on the $8$-letter alphabet $\{s, n, o, w, s', n', o', w'\}$
by 
$$
\begin{array}{llllll}
s  &\to &s  \ &n  \ &o  \ &w  \\
w  &\to &s  \ &n  \ &o  \ &w  \\
o  &\to &o  \ &w  \ &n' \ &o  \\
n  &\to &n  \ &o' \ &s  \ &n  \\
s' &\to &s' \ &n' \ &o' \ &w' \\
w' &\to &s' \ &n' \ &o' \ &w' \\ 
o' &\to &o' \ &w' \ &n  \ &o' \\
n' &\to &n' \ &o  \ &s' \ &n' \\
\end{array}
$$
followed by the map $s \to 1$, $w \to 1$, $n \to -j^2$, $o \to -j$,
$s' \to -1$, $w' \to -1$, $n' \to j^2$, $o' \to j$ (actually there is a 
misprint in Dekking's paper for the images of $o$ and $o'$ by this map).
Dekking's morphism clearly gives the same fixed point beginning with $s$
as the morphism
$$
\begin{array}{llllll}
s  &\to &s  \ &n  \ &o  \ &s  \\
o  &\to &o  \ &s  \ &n' \ &o  \\
n  &\to &n  \ &o' \ &s  \ &n  \\
s' &\to &s' \ &n' \ &o' \ &s' \\
o' &\to &o' \ &s' \ &n  \ &o' \\
n' &\to& n' \ &o  \ &s' \ &n' \\
\end{array}
$$
followed by the map $s \to 1$, $n \to -j^2$, $o \to -j$, $s' \to -1$,
$n' \to j^2$, $o' \to j$. This last map being a bijection from
$\{s, n, o, s', n', o'\}$ to $\{1, -1, j, -j, j^2, -j^2\}$, we can write
this morphism directly on the alphabet $\{1, -1, j, -j, j^2, -j^2\}$ without
changing the final sequence:
$$
\begin{array}{rrrrrr}
1    &\to &1    \ &-j^2 \ &-j   \ &1    \\
-j   &\to &-j   \ &1    \ &-j^2 \ &-j   \\
-j^2 &\to &-j^2 \ &j    \ &1    \ &-j^2 \\
-1   &\to &-1   \ &j^2  \ &j    \ &-1   \\
j    &\to & j   \ &-1   \ &-j^2 \ &j    \\
j^2  &\to & j^2 \ &-j   \ &-1   \ &j^2  \\
\end{array}
$$
(this morphism is the morphism $\lambda$ we obtained above up to changing
$j$ into $j^2$).

\end{itemize}

\section{A generalization}

Let $p$ and $q$ be two integers $\geq 2$. Define $s_p(n)$ to be
the sum of digits in the $p$-ary expansion of the integer $n$.
Define $\eta_p := e^{2i\pi/p}$ and $\eta_q = e^{2i\pi/q}$. 
The general sum
$$
Z(N, p, q) := \sum_{0 \leq k < N} \eta_p^{s_p(k)} \eta_q^k
$$
was considered by Dekking in \cite{Dekking-bx}.
In the case where $p=2$, hence $\eta_p = -1$ and $s_p = s_2$, the 
same procedure as in \cite{Coquet, Dekking1} can be applied, replacing 
the above morphism of length $4$
$$
z \ \to \ z \ \ -j z \ \ -j^2 z \ \ \ z
$$
by the morphism of length $Q$
$$
z \ \ \to \ \ z \ \ \ \ \ (-1)^{s_2(1)}\eta_q z \ \ \ \ \ 
(-1)^{s_2(2)}\eta_q^2 z \ \ \ \ \ \cdots \ \ \ \ \ 
(-1)^{s_2(Q-1)}\eta_q^{Q-1}z 
$$
where $Q := 2^{\varphi(q)}$ and $\varphi$ is the Euler function.
All details are given in \cite{Dekking-bx}, where Dekking gives in 
particular (see \cite[Fig.~1, p.~32-06]{Dekking-bx}) fractal pictures
associated with $p=2$ and $q \in \{3, 5, 7, 13\}$: the case $q=3$,
i.e., $Q = 4$ and $\eta_q = j$, shows of course the von Koch curve.

\section{Automatic sequences and fractal objects}

The purpose of this section is to emphasize that there is a general link
between automatic sequences and fractal objects. A very first remark is 
that being a fixed point of a (uniform) morphism means having some sort 
of self-similarity: since an automatic sequence is infinite, a rescaling 
process is often necessary to ``see'' the fractal object. Let us begin with 
a very simple example. Define the morphism of length $3$ on the alphabet 
$\{0, 1\}$ by
$$
0 \to 000 \ \ \ \ \ 1 \to 101.
$$
Its fixed point beginning in $1$ is obtained as the limit of the
following sequence of words
$$
\begin{array}{ll}
&1 \\
&1 \ 0 \ 1 \\
&1 \ 0 \ 1 \ 0 \ 0 \ 0 \ 1 \ 0 \ 1 \\
&1 \ 0 \ 1 \ 0 \ 0 \ 0 \ 1 \ 0 \ 1 \ 0 \ 0 \ 0 \ 0 \ 0 \ 0 \ 0 \ 0 \ 0 \
1 \ 0 \ 1 \ 0 \ 0 \ 0 \ 1 \ 0 \ 1 \\
\end{array}
$$
The reader will have immediately noted that, starting with $101$, representing 
$1$'s by black horizontal segments of length $1/3$ and $0$'s by white horizontal 
segments of length $1/3$, and ``renormalizing'' by a factor of $1/3$ at each step
yields the classical ternary Cantor set. This approach can be generalized
to situations where the representation does not necessarily lie on a straight
line: see in particular the papers of Golomb \cite{Golomb} and Giles
\cite{Giles1, Giles2, Giles3} where fixed points of morphisms and automatic 
sequences are not mentioned but can be found, and the papers of Dekking 
\cite{Dekking1, Dekking2, Dekking-bx}. In the same direction see the papers
of Dekking, Dekking and Mend\`es France, Deshouillers, Mend\`es France, and 
Mend\`es France and Shallit \cite{Dekking3, DMF, Deshouillers, Mendes, MFS}, 
but also the paper of Salon with an appendix by Shallit \cite{Salon}, where 
two-dimensional automatic sequences enter the picture, without forgetting the 
article by Prusinkiewicz and Hammel \cite{PH}, and the article by Allouche, 
Allouche and Shallit \cite{AAS}. We would also like to quote two papers by
von Haeseler, Peitgen, Skordev \cite{HPS1, HPS2} and the references therein,
and a paper by Barb\'e and von Haeseler \cite{BH}.

\bigskip

Another relation between fractal objects and fixed points of morphisms or 
automatic sequences $(a_n)_{n \geq 0}$ occurs in the study of the summatory 
function $\sum_{n \leq x} a_n$, whose asymptotic behavior often involves 
a logarithmically periodic, continuous and nowhere (or almost nowhere) 
differentiable function: this non-regular behavior corresponds to
a geometric object which is fractal. A nice and usually difficult question
is to determine the Holder exponents of the non-regular terms in the
asymptotic behavior of these summatory functions. While the first such 
study, due to Delange \cite{Delange}, involves the summatory function of 
the sequence $(s_2(n))_{n \geq 0}$, several papers deal with the summatory 
function of fixed points of morphisms or automatic sequences, or even with 
more general classes of sequences: see in particular the work of Dumont, 
and Dumont and Thomas, in particular \cite{Dumont, DT}, the thesis of 
Cateland \cite{Cateland}, the paper by Flajolet, Grabner, Kirschenhofer, 
Prodinger, and Tichy \cite{FGKPT}, and finally the paper of Tenenbaum 
\cite{Tenenbaum}.

\bigskip

A third direction we want to mention concerns noiselets with in particular
an occurrence of the ($2$-automatic) Shapiro-Rudin sequence, see \cite{CGM}.

\bigskip

What precedes should not hide the specificity of the case addressed here,
namely the close relation between the von Koch curve and the Thue-Morse
sequence. While the relations explained in this section are quite general,
they either link a well-known fractal object to an {\it ad hoc\,} automatic
sequence, or a well-known automatic sequence gives rise to some unnamed
and not ``classical'' fractal object. In the case of von Koch and Thue-Morse,
the fractal {\it and\,} the sequence were studied quite independently, before
a link was found more then eighty years after von Koch's paper and more than
seventy years after Thue's paper.

\bigskip

\noindent
{\bf Acknowledgments} This paper was written during two stays at Cevis-Mevis
of JPA who thanks heartily all the colleagues in Bremen, and in particular 
H.-O. Peitgen, for their very friendly hospitality. Both authors thank the
referee for pertinent and useful comments.

\end{document}